\documentclass{amsart}
\begin{document}

\newcommand{\ea}{\mbox{{\bf a}}}
\newcommand{\eu}{\mbox{{\bf u}}}
\newcommand{\ueu}{\underline{\eu}}
\newcommand{\ueo}{\overline{u}}
\newcommand{\oeu}{\overline{\eu}}
\newcommand{\ew}{\mbox{{\bf w}}}
\newcommand{\ef}{\mbox{{\bf f}}}
\newcommand{\eF}{\mbox{{\bf F}}}
\newcommand{\eC}{\mbox{{\bf C}}}
\newcommand{\en}{\mbox{{\bf n}}}
\newcommand{\eT}{\mbox{{\bf T}}}
\newcommand{\eL}{\mbox{{\bf L}}}
\newcommand{\eR}{\mbox{{\bf R}}}
\newcommand{\eV}{\mbox{{\bf V}}}
\newcommand{\eU}{\mbox{{\bf U}}}
\newcommand{\ev}{\mbox{{\bf v}}}
\newcommand{\eve}{\mbox{{\bf e}}}
\newcommand{\uev}{\underline{\ev}}
\newcommand{\eY}{\mbox{{\bf Y}}}
\newcommand{\eK}{\mbox{{\bf K}}}
\newcommand{\eP}{\mbox{{\bf P}}}
\newcommand{\eS}{\mbox{{\bf S}}}
\newcommand{\eJ}{\mbox{{\bf J}}}
\newcommand{\eB}{\mbox{{\bf B}}}
\newcommand{\eH}{\mbox{{\bf H}}}
\newcommand{\leb}{\mathcal{ L}^{n}}
\newcommand{\eI}{\mathcal{ I}}
\newcommand{\eE}{\mathcal{ E}}
\newcommand{\hen}{\mathcal{H}^{n-1}}
\newcommand{\eBV}{\mbox{{\bf BV}}}
\newcommand{\eA}{\mbox{{\bf A}}}
\newcommand{\eSBV}{\mbox{{\bf SBV}}}
\newcommand{\eBD}{\mbox{{\bf BD}}}
\newcommand{\eSBD}{\mbox{{\bf SBD}}}
\newcommand{\ecs}{\mbox{{\bf X}}}
\newcommand{\eg}{\mbox{{\bf g}}}
\newcommand{\paromega}{\partial \Omega}
\newcommand{\gau}{\Gamma_{u}}
\newcommand{\gaf}{\Gamma_{f}}
\newcommand{\sig}{{\bf \sigma}}
\newcommand{\gac}{\Gamma_{\mbox{{\bf c}}}}
\newcommand{\deu}{\dot{\eu}}
\newcommand{\dueu}{\underline{\deu}}
\newcommand{\dev}{\dot{\ev}}
\newcommand{\duev}{\underline{\dev}}
\newcommand{\weak}{\rightharpoonup}
\newcommand{\weakdown}{\rightharpoondown}
\renewcommand{\contentsname}{ }

\newcommand{\down}[1]{{#1}^{\downarrow}}
\newcommand{\up}[1]{{#1}^{\uparrow}}
\newcommand{\fromto}[2]{\left\{{#1},...,{#2}\right\}}
\newcommand{\lwd}{\prec_{w}}
\newcommand{\lwu}{\prec^{w}}

\newtheorem{rema}{Remark}[section]
\newtheorem{thm}{Theorem}[section]
\newtheorem{lema}{Lemma}[section]
\newtheorem{prop}{Proposition}[section]
\newtheorem{defi}{Definition}[section]
\newtheorem{conje}{Conjecture}[]
\newtheorem{exempl}{Example}[section]
\newtheorem{opp}{Open Problem}[]
\renewcommand{\contentsname}{ }
\newenvironment{rk}{\begin{rema}  \em}{\end{rema}}
\newenvironment{exemplu}{\begin{exempl}  \em}{\end{exempl}}

\title{Majorisation with applications to the Calculus of 
Variations}
\author{MARIUS BULIGA}

%\markboth{Marius Buliga}{Majorisation and Calculus of Variations}

\date{\today}

\begin{abstract}
This paper explores some connections between rank one convexity, 
multiplicative quasiconvexity and Schur convexity. Theorem 
\ref{trank1} gives simple necessary and sufficient conditions for 
an isotropic objective function to be rank one convex on the set 
of matrices with positive determinant. Theorem \ref{texe} describes a 
class of possible non-polyconvex but multiplicative quasiconvex 
isotropic functions. This class is not contained in a well known 
theorem of Ball (\ref{critball} in this paper) which gives sufficient 
conditions for an isotropic and objective function 
to be polyconvex. We show here that there is a new way to prove 
directly the quasiconvexity (in the multiplicative form). Relevance 
of Schur convexity for the description of rank one convex hulls is 
explained. 
\end{abstract}


\maketitle


\section{Introduction}

There is a strong resemblance between the following two theorems. 
The first theorem is (Horn, \cite{horn}(1954), Thompson 
\cite{thompson}(1971), theorem 1.):

\begin{thm}
Let $X, Y$ be any two positive definite $n \times n$ matrices and let 
$x_{1} \geq x_{2} \geq ... \geq x_{n}$ and 
$y_{1} \geq y_{2} \geq ... \geq y_{n}$ 
denote the respective sets of eigenvalues. Then there is an unitary 
matrix $U$ such that $XU$ and $Y$ have the same spectrum if and only 
if: 
$$\prod_{i=1}^{k} x_{i} \geq \prod_{i=1}^{k} y_{i} \ \ , \ k = 
1, ... , n-1$$
$$\prod_{i=1}^{n} x_{i} \ = \ \prod_{i=1}^{n} y_{i}$$
\label{hornthom}
\end{thm}

The second theorem is (Dacorogna, Tanteri  \cite{dacotan}(2001), 
theorem 20, see also Dacorogna, Marcellini \cite{dacomarc}):

\begin{thm}
Let $0 \leq \lambda_{1}(A) \leq ... \leq \lambda_{n}(A)$ denote the 
singular values of a matrix $A \in R^{n \times n}$ and 
$$E(a) = \left\{ A \in R^{n\times n} \ : \ \lambda_{i}(A) = a_{i}
 \ , \ 
i = 1,...,n \ , \  \det A = \prod_{i=1}^{n} a_{i}\ \right\}$$
The following then holds
$$Pco \ E =   \ Rco \ E(a) = \left\{ A \in R^{n \times n} \ : \ 
\prod_{i=\nu}^{n} \lambda_{i}(A) \leq \prod_{i=\nu}^{n} a_{i} \ , \ 
\nu = 2, ... ,n \ , \  \det A = \prod_{i=1}^{n} a_{i}  \ \right\}$$
where $PCo$, $Rco$ stand for polyconvex, rank one convex envelope. 
\label{dacotan}
\end{thm}

Both theorems can be understood as describing the set 
$\left\{ y \ : \ y\prec \prec x \right\}$ where $\prec \prec$ is 
some order relation connected to the products which rise in each 
theorem. 

It turns out that a common frame of these apparently 
scattered results is the notion of majorisation. Majorisation comes 
in pair with Schur convexity. The purpose of this note is to study 
the monotonicity (Schur convexity in particular) properties 
of rank one convex functions with respect to majorisation relation. 

The content of the paper is described further. After the setting of 
notations in section 2, section 3 gives a brief passage trough basic
 properties of the majorisation relations. Section 4 lists some 
properties of singular values and eigenvalues of matrices connected 
to majorisation. Section 5 concerns simple necessary and 
sufficient conditions for an objective, isotropic function to be 
rank one convex on the set of matrices with positive determinant. 
In section 6 is described a class 
of objective isotropic functions which are multiplicative 
quasiconvex. This class is interesting because most of it's 
elements seem to be non-polyconvex. This class is not contained in 
a well known theorem of Ball (\ref{critball} in this paper) which 
gives sufficient conditions for an isotropic and objective function 
to be polyconvex. It turns out that 
Schur convexity can be used to prove quasiconvexity.  In section 7 
the resemblance between theorems \ref{hornthom} and \ref{dacotan} 
is explained.

\section{Notations}

\begin{tabular}{ll}
$A$,$B$, ... & real or complex matrices \\ 
$x$,$y$,$u$,$v$, ... & real or complex vectors \\ 
$\lambda(A)$ & the vector of eigenvalues of $A$ \\
$\sigma(A)$ & the vector of singular values of $A$ \\
$A^{*}$ & the conjugate transpose of $A$ \\
$A^{T}$ & the transpose of $A$ \\
$diag(A)$ & the diagonal of $A$, seen as a vector \\
$Diag(v)$ & the diagonal matrix constructed from the vector $v$ \\
$S_{n}$ & the set of permutation (of coordinates) matrices\\
$conv(A)$ & the convex hull of the set $A$\\
$\circ$ & function composition
\end{tabular}

For any matrix $A \in gl(n,C)$, the matrix $A^{*}A$ is Hermitian. The 
eigenvalues of the square root of $A^{*}A$ are, by definition, the 
singular values of $A$.

Matrices are identified with linear transformations.

For a vector $x \in R^{n}$ we denote by $\down{x}$, $\up{x}$,  the 
vectors obtained by rearranging the coordinates of $x$ in 
decreasing, respectively increasing orders.

\section{Basics about majorisation}

We have used Bhatia \cite{bha}, Chapter 2, and Marshall, Olkin 
\cite{mo}, Chapters 1-3. The results are given in the logical order. 

\begin{defi}
The following majorisation notions are partial  order relations  in 
$R^{n}$.Let  $x,y \in R^{n}$ be arbitrary vectors. Then:
\begin{enumerate}
\item[$\bullet$] $x \leq y$ if $x_{i} \leq y_{i}$ for any $i \in 
\fromto{1}{n}$.
\item[$\bullet$] $x \lwd y$ if 
$$ \sum_{j=1}^{k} \down{x_{j}} \leq \sum_{j=1}^{k} \down{y_{j}}$$
for any $k \in \fromto{1}{n}$. We say that $x$ is submajorised by 
$y$.  
\item[$\bullet$] $x\prec y$ if $x \lwd y$ and 
$$ \sum_{j=1}^{n} \down{x_{j}} = \sum_{j=1}^{n} \down{y_{j}}$$
We say that $x$ is majorised by $y$. 
\end{enumerate}
\label{d1}
\end{defi}

The notion of majorisation, the last in definition \ref{d1}, is the 
most interesting. See Marshall, Olkin \cite{mo}, Chapter 1, for the 
various places when one can encounter it. The majorisation is  in 
closed relationship with the notions 
of a T-transform and a doubly-stochastic matrix. 

\begin{defi}
A linear map $A$ on $R^{n}$ is a T-transform if there are $t \in 
[0,1]$, 
$i,j \in \fromto{1,}{n}$ such that $$(Ax)_{k} = x_{k}$$ for any 
$k$ different from $i,j$, $$(Ax)_{i} = tx_{i}+(1-t)x_{j} \ , \ 
(Ax)_{j} = tx_{j}+(1-t)x_{i}$$
A  matrix $A \in gl(n,R)$ is called doubly-stochastic if $$A_{ij}
\geq 0$$ 
$$\sum_{k=1}^{n}A_{kj} = 1 \ , \ \sum_{k=1}^{n}A_{ik} = 1$$ 
for all $i,j$. Any T-transform is doubly-stochastic. Matrices which 
correspond to permutation of coordinates are also doubly-stochastic. 
\label{d2}
\end{defi} 

The property of a matrix of being doubly stochastic can be 
formulated in terms of majorisation. 

\begin{thm} 
A matrix $A$ is doubly stochastic if and only if $Ax \prec x$ for 
any $x \in R^{n}$. 
\label{t1}
\end{thm}

Conversely, we have: 

\begin{thm}
The following statements are equivalent: 
\begin{enumerate}
\item[(i)] $x \prec y$
\item[(ii)] $x$ is obtained from $y$ by a finite number of 
T-transforms
\item[(iii)] $x = Ay$ for some doubly stochastic matrix $A$
\end{enumerate}
\label{t2}
\end{thm}

There are strong connections between majorisation, doubly 
stochastic matrices 
and convexity. These will make the subject relevant for the 
Calculus of Variations further.

\begin{thm}
(Birkhoff) The set of doubly stochastic matrices is the convex 
hull of the set of permutation matrices. 
\label{t3}
\end{thm}

\begin{thm}(Hardy, Littlewood, Polya) 
The following statements are equivalent: 
\begin{enumerate}
\item[(i)] $x \prec y$
\item[(ii)] $x$ is in the convex hull of $S_{n} x$
\item[(iii)] for any convex function 
$\phi$ from $R$ to $R$ we have 
$$\sum_{i=1}^{n} \phi(x_{i}) \leq \sum_{i=1}^{n} \phi(y_{i})$$
\end{enumerate}
\label{t4}
\end{thm}

With any order relation comes an associated monotonicity notion. 

\begin{defi}
Consider a map $\Phi$ defined from an $S_{n}$ invariant set in 
$R^{n}$, with range in $R^{m}$. We say that $\Phi$ is: 
\begin{enumerate}
\item[$\bullet$]  increasing if 
$$x\leq y \ \Longrightarrow \Phi(x) \leq \Phi(y)$$
\item[$\bullet$] convex if for all $t \in [0,1]$ 
$$\Phi(tx+(1-t)y) \leq t\Phi(x)+(1-t)\Phi(y)$$
\item[$\bullet$] isotone if 
$$x\prec y \ \Longrightarrow \Phi(x) \lwd \Phi(y)$$
\item[$\bullet$] strongly isotone if
$$x\lwd y \ \Longrightarrow \Phi(x) \lwd \Phi(y)$$
\item[$\bullet$] strictly isotone if 
$$x\prec y \ \Longrightarrow \Phi(x) \prec \Phi(y)$$
\end{enumerate}
Any isotone $\Phi$ with range in $R$ is called Schur-convex. Note 
that convexity in the sense of this definition matches with the 
classical notion for functions $\Phi$ with range in $R$. 
\label{d3}
\end{defi}

The next theorem shows that symmetric convex maps are isotone. 

\begin{thm}
Let $\Phi:R^{n} \rightarrow R^{m}$ be convex. If for any $P \in 
S_{n}$ there is $P' \in S_{m}$ such that 
$\Phi \circ P = P'\circ \Phi$ then $\Phi$ is isotone. 
If in addition $\Phi$ is monotone increasing then $\Phi$ is 
strictly isotone. 
\label{t5}
\end{thm}
In particular any $L^{p}$ norm on $R^{n}$ is Schur-convex. Not any 
isotone function is convex, though. Important examples are the 
elementary symmetric polynomials, which are not convex but they are 
{\it Schur-concave}. 

One can give three characterizations of isotone (or Schur convex) 
functions $f:R^{n} \rightarrow R$. Before that we need some 
notations. 

Let us begin by noticing  that the permutation group 
$S_{n}$ acts on $GL(n,R)^{+}$ in this obvious way: for any 
$P \in S_{n}$ and any $F \in GL(n,R)^{+}$ the matrix 
$P.F \in GL(n,R)^{+}$ has components  
$(P.F)_{ij} = F_{P(i)P(j)}$.

Let 
$$\mathcal{D} = \left\{ x \in R^{n} \ : \ x_{1} \geq x_{2} 
\geq ...\geq x_{n} 
\right\}$$
We shall call a function $f:A \subset R^{n} \rightarrow R$ 
symmetric if for any 
permutation matrix $P \in S_{n}$ $P(A) \subset A$ and 
$f\circ P = f$. The partial derivative of $f$ with respect to 
$x_{i}$ will be denoted by $f_{i}$.

\begin{thm}
Let $I$ be an open interval in $R$ and let 
$f:I^{n} \rightarrow R$ be continuously differentiable. Then 
$f$ is Schur convex if and only if one of the following conditions 
is true: 
\begin{enumerate}
\item[(a)] (Schur)  $f$ is symmetric an $f_{i}$ is decreasing in 
$x_{i}$ for all $x \in \mathcal{D} \cap I^{n}$. 
\item[(b)] (Schur) $f$ is symmetric and for all $i \not = j$ 
$$(x_{i} -x_{j})(f_{i}(x) - f_{j}(x)) \geq 0$$
\end{enumerate}
Eliminate from the hypothesis the differentiability of $f$ and 
consider 
$$f: A \subset R^{n} \rightarrow R$$ with $A$ symmetric. Then $f$ is 
Schur convex if and only if: 
\begin{enumerate}
\item[(c)] $f$ is symmetric and  
$$x_{1} \mapsto f(x_{1}, s - x_{1}, x_{2}, ... , x_{n})$$
is increasing in $x_{1} \geq s/2$, for any fixed $s,x_{3}, .... , 
x_{n}$. 
\end{enumerate}
\label{t1schur}
\end{thm}

For weak majorisation and strongly isotone functions 
we have the following theorem: 

\begin{thm}
Let $I$ be an open interval in $R$ and let 
$f:I^{n} \rightarrow R$.
\begin{enumerate}
\item[(a)] (Ostrowski)Let $f$ be continuously differentiable. Then 
$f$ is strongly isotone  if and only if    $f$ is symmetric and  
for all $x \in \mathcal{D} \cap I^{n}$ we have $Df(x) \in 
\mathcal{D} \cap R^{n}_{+}$, that is: 
$$f_{,1}(x) \geq f_{,2}(x) \geq ... \geq f_{,n}(x) \geq 0$$
\item[(b)] Without differentiability assumptions, $f$ is strongly 
isotone if and only if $f$ increasing and Schur convex. 
\end{enumerate}
\label{tostrow}
\end{thm}

\section{Order relations for matrices}

The results from this section have deep connections with  Lie group 
theory, symplectic geometry and sub-Riemannian geometry. I shall 
give here only a minimal presentation, for  matrix groups.  

Main references are again Bhatia \cite{bha}, Chapter 2, and 
Marshall, Olkin \cite{mo}, Chapter 3; also Thompson 
\cite{thompson}. The paper Kostant \cite{konstant} gives an image 
of what's really happening from the Lie group point of view. 

\begin{defi}
We denote by $\mathcal{P}(n)$ the cone of Hermitian, positive 
definite matrices. 
In the class of Hermitian matrices we have the order relation 
$A \geq B$ if $A-B \in \mathcal{P}(n)$. 
\end{defi}

The order relation $\leq$ between Hermitian matrices reflects into 
the order relation between the eigenvalues. The next theorem, 
belonging to Weyl, is theorem F1, chapter 16, Marshall, Olkin 
\cite{mo}. 

\begin{thm}
(Weyl) If $A,B$ are Hermitian matrices such that $A \leq B$ then
$$\down{\lambda}(A) \leq \down{\lambda}(B)$$
\label{tf1}
\end{thm}

For $A,B$ matrices, their Schur product is the matrix $A \odot B$ 
given by 
$$(A \odot B)_{ij} = A_{ij} B_{ij} \mbox{ (no summation)}$$

\begin{thm}
(Schur) If $A \geq B$ and $C\geq 0$ then $A\odot C \geq B \odot C$. 
\end{thm}

Next theorem shows a first connection between majorisation and 
symmetric matrices. 

\begin{thm}
(Schur) For any symmetric matrix $A$ we have $diag(A) \prec 
\lambda(A)$. 

(Horn) Conversely, given vectors $a, b \in R^{n}$ such that 
$a \prec b$, there is 
a symmetric matrix $A$ such that $diag(A) = a$ and $\lambda(A) = b$
\label{tschurhorn}
\end{thm}

In $GL(n,R)^{+}$ we have the order relation (introduced by 
Thompson \cite{thompson})
$$X \prec Y \ \ \mbox{ if } \log \sigma(X) \prec \log \sigma(y)$$

Horn-Thompson theorem \ref{hornthom} can be reformulated as: 

\begin{thm}
(Horn,Thompson, theorem \ref{hornthom} reformulated) 
Let $X, Y$ be any two positive definite $n \times n$ matrices and 
let 
$x_{1} \geq x_{2} \geq ... \geq x_{n}$ and $y_{1} \geq y_{2} 
\geq ... \geq y_{n}$ 
denote the respective sets of eigenvalues. Then there is an 
unitary matrix $U$ such that $XU$ and $Y$ have the same spectrum 
if and only if  $Y \prec X$. 
\end{thm}

Another interesting majorisation occurs between the absolute value 
of eigenvalues and singular values respectively. 

\begin{thm}
(Weyl) For any matrix $F \in GL(n,C)$ we have the inequality: 
$$\log \mid \lambda(F) \mid \ \prec \ \log \sigma(F)$$
\label{tweyl}
\end{thm}

Finally, denote by $\mid A \mid$ the spectral radius of $A$, i.e. 
the maximum over the modulus of singular values of $A$. Next 
theorem shows the algebraic deep of Thompson's order relation. 

\begin{thm}
(Kostant) $X \prec Y$ if and only if for any linear representation 
$\pi$ of the group $GL(n,R)^{+}$ we have $\mid \pi(X) \mid \leq 
\mid \pi(Y) \mid$. 
\label{tkons}
\end{thm}

\section{Objective isotropic elastic potentials}
We are interested in functions $w: GL(n,R)^{+} \rightarrow R$  
which are objective 
$$\forall \ Q \in SO(n) \ w(QF) = W(F)$$
and isotropic 
$$\forall \ Q \in SO(n) \ w(FQ) = W(F)$$

If $w$ is $C^{2}$, then we call it rank one convex if it satisfies 
the ellipticity condition: 
\begin{equation}
\sum_{i,j,k,l = 1}^{n} \frac{\partial^{2}w}{\partial F_{ij} 
\partial F_{kl}} (F) 
a_{i}b_{j}a_{k}b_{l} \ \geq \ 0
\label{elipc}
\end{equation}
for any $F \in GL(n,R)^{+}$, $a,b \in R^{n}$. 

For $A,B$ matrices, we denote by $[[A,B]]$ the segment 
$$[[A,B]] \ = \ \left\{ (1-t)A + tB \ : \ t \in [0,1] \right\}$$

We have the more general definition of rank one convexity: 
\begin{defi}
The function $w: GL(n,R)^{+} \rightarrow R$ is rank one convex if 
for any $A,B \in GL(n,R)^{+}$ such that $rank \ (A-B) = 1$ and 
$[[A,B]] \subset GL(n,R)^{+}$ the function $t \in [0,1] \mapsto 
w((1-t)A + tB) \in R$ is convex. 
\label{drk1}
\end{defi}

It is straightforward that the ellipticity condition is equivalent 
to rank one convexity for smooth functions. 

If $w$ is objective and isotropic then there is a symmetric function 
$g: R^{n}_{+} \rightarrow R$ such that $w(F) = g(\sigma(F))$. If 
$w$ is $C^{2}$ then $g$ is too. 

We shall introduce two auxiliary functions. They are the following: 
$$h: R^{n} \rightarrow R \ \ , \ \ h(x) = g(\exp x)$$
$$l: R^{n}_{+} \rightarrow R \ \ , \ \ l(x) = g(\sqrt{x})$$
The auxiliary function $h$ will be called "the diagonal of $w$".

Now, let $x \in R^{n}$ or $R^{n}_{+}$ such that $x_{i} \not = 
x_{j}$ for $i \not = j$. 
The following quantities will help. 
$$\Gamma_{ij}(x) = \frac{h_{i}(x) - h_{j}(x)}{x_{i} - x_{j}}$$
Any function $\Gamma_{ij}$ can be prolonged by continuity in 
$x_{i} = x_{j}$. In order to shorten the notation we shall put in 
the functions arguments only the 
terms that count. For example $f(x_{i}, x_{j})$ means $f(x)$ and 
$f(x_{i}, x_{i})$ is $f(x)$ for an $x$ such that $x_{i} = x_{j}$. 
With this notation one can define by continuity: 
$$\Gamma_{ij}(x_{i}, x_{i}) = h_{ij}(x_{i}, x_{i}) - h_{jj}(x_{i}, 
x_{i})$$
A straightforward computation shows that 
$$\Gamma_{ij}(x_{i},x_{j}) = G_{ij}(\exp x_{i}, \exp x_{j})
 (\exp x_{i} + \exp x_{j})$$
where 
$$G_{ij}(x_{i}, x_{j}) = \frac{x_{i} g_{i}(x_{i},x_{j}) - x_{j} 
g_{j}(x_{i},x_{j})}{x_{i}^{2} - x_{j}^{2}}$$
Therefore the coefficients $G_{ij}$ can be prolonged by continuity 
to $x_{i} = x_{j}$. We shall put $\Gamma_{ii} = G_{ii} = 0$. 

We shall introduce also the symmetric matrix $\Xi$. We shall define 
it first for vectors $x \in R^{n}_{+}$ with all components 
different and then extend it by continuity to all vectors. 
For $i \not = j$ we define: 
$$\Xi_{ij}(x_{i}, x_{j}) = \frac{l_{i}(x_{i}, x_{j}) - l_{j}
(x_{i},x_{j})}{x_{i} - x_{j}}$$
If $x_{i} = x_{j}$ then the prolongation by continuity of 
$\Xi_{i,j}$ is 
$$\Xi_{ij}(x_{i}, x_{i}) = l_{ii}(x_{i},x_{i}) - l_{ij}(x_{i},
x_{i})$$
For $i=j$ and all $x$, we define: 
$$\Xi_{ii}(x) = 2 l_{ii}(x) + l_{i}(x)/x_{i}$$
Again by straightforward computation we find that, for any $x$ 
with all components different: 
$$\Xi(x) \odot \left( \sqrt{x} \otimes \sqrt{x} \right) \ = \ 
\bar{H}(\sqrt{x})$$
where the matrix $\overline{H}$ is defined by: 
$$\overline{H}_{ij}(x_{i}, x_{j}) = \frac{x_{j}g_{i}(x_{i},x_{j}) - 
x_{i}g_{j}(x_{i},x_{j})}{x_{i}^{2} - x_{j}^{2}}$$
for $i \not = j$ and 
$$\overline{H}_{ii}(x) = 0$$
As previously, the function $\overline{H}_{ij}$ can be prolongated 
by continuity to $x_{i} = x_{j}$.

A consequence of theorem 6.4 Ball \cite{ball1} is: 

\begin{prop}
 For $x$ with all components different, 
the ellipticity condition \eqref{elipc} for the objective isotropic 
function $w$  can be expressed in terms of the associated function 
$g$ as  
$$\sum_{i,j=1}^{n} g_{ij} a_{i} a_{j} b_{i} b_{j} \ + \ \sum_{i 
\not = j} 
G_{ij}a_{i}^{2}b_{j}^{2} \ + \ \sum_{i \not = j} \overline{H}_{ij} 
a_{i} a_{j} b_{i} b_{j} \ \geq 0$$
\label{pball}
\end{prop}
By continuity arguments it follows that one can write the
 ellipticity condition for all $x \in R^{n}_{+}$ like this:
\begin{equation} 
\sum_{i,j=1}^{n} H_{ij} a_{i} a_{j} b_{i} b_{j} \ + \ 
\sum_{i,j = 1}^{n} 
G_{ij}a_{i}^{2}b_{j}^{2} \ \geq 0
\label{mainc}
\end{equation}
where $H$ is the matrix $H= \overline{H} + D^{2}g$.

\begin{thm}
Necessary and sufficient conditions for $w \in C^{2}$ to be rank 
one convex are: 
\begin{enumerate}
\item[(a)] $h$ is Schur convex and 
\item[(b)] for any $x \in R^{n}$ we have 
\begin{equation}
H_{ij} x_{i}x_{j} + G_{ij}\mid x_{i}\mid \mid x_{j} \mid \ \geq 0
\label{mainineq}
\end{equation} 
\end{enumerate}
\label{trank1}
\end{thm}

\begin{rk}
The condition (a) is equivalent to the Baker-Ericksen \cite{be} 
set of inequalities
$$\frac{x_{i}g_{i}(x_{i},x_{j}) - x_{j}g_{j}(x_{i},x_{j})}{x_{i}^{2}
 - x_{j}^{2}} \geq 0$$
for all $i\not= j$ and $x_{i} \not = x_{j}$. Indeed, by theorem 
\ref{t1schur} (b), the function $h$ is Schur convex if and only if 
$$\left(h_{i}(x_{i},x_{j}) - h_{j}(x_{i},x_{j})\right)\left(x_{i} - 
x_{j}\right) \geq 0$$
for all $i\not= j$ and $x_{i} \not = x_{j}$. But the definition of 
the "diagonal" $h$ and obvious computation show the equivalence 
between the two sets of inequalities.
Silhavy \cite{silhavy} expresses Baker-Ericksen inequalities using 
multiplication instead division, but apparently he does not make 
this obvious connection with Schur convexity.  
\end{rk}

\begin{proof}
We prove first the sufficiency. The hypothesis is that for all 
$i,j$ $G_{ij}\geq 0$ 
and for all $x \in R^{n}$ the relation \eqref{mainineq} holds. We 
claim that for any $a,b \in R^{n}$ the inequality 
$$G_{ij}a_{i}a_{j}b_{i}b_{j} \leq G_{ij}a_{i}^{2}b_{j}^{2}$$
is true. The ellipticity condition follows then from 
\eqref{mainineq} by the choice 
$x_{i}=a_{i}b_{i}$, for each $i = 1,...,n$. Indeed, we have the 
chain of inequalities
$$0\leq H_{ij}a_{i}b_{i}a_{j}b_{j} + G_{ij}\mid a_{i}b_{i}
\mid \mid a_{j}b_{j}\mid \ \leq \ H_{ij}a_{i}b_{i}a_{j}b_{j} + 
G_{ij} a_{i}^{2}b_{j}^{2}$$
In order to prove the claim note that $G_{ij}\geq 0$ implies
$$-G_{ij}(a_{j}b_{i} - a_{i}b_{j})^{2} \ \leq 0$$
A straightforward computation which uses the relations $G_{ij}=
G_{ji}$ gives 
$$0 \geq \ -G_{ij}(a_{j}b_{i} - a_{i}b_{j})^{2} \ = \ 
2 G_{ij} (a_{j}b_{i} - a_{i}b_{j})a_{i}b_{j}$$
The sufficiency part is therefore proven. 

For the necessity part choose first in the ellipticity condition 
$a_{i} = \delta_{iI}$, $b_{i} = \delta_{iJ}$. For $I\not=J$ we obtain 
$G_{ij} \geq 0$, which means the Schur convexity of $h$. (For $I=J$ 
we obtain $g_{ii} \geq 0$, interesting but with no use in this 
proof.) 

Next, suppose that $x, a \in (R^{*})^{n}$ and choose $b_{i} = 
x_{i}/a_{i}$ for each $i = 1,...,n$. The ellipticity condition gives: 
$$\sum_{i,j} H_{ij}x_{i}x_{j} \ + \ \sum_{i,j=1}^{n}G_{ij}
\left(\frac{a_{i}}{a_{j}}\right)^{2}x_{j}^{2} \ \geq 0$$
Take $a_{i}^{2} = \mid x_{i} \mid$ and get \eqref{mainineq}, but 
only for $x \in (R^{*})^{n}$. The expression from the left of 
\eqref{mainineq} makes sense 
for any $x$. By continuity with respect to $x$ we prove the thesis. 
\end{proof}

There is a certain interest in giving necessary and sufficient 
conditions for an objective isotropic $w$ to be rank one convex, 
especially in the cases $n=2$ and $n=3$. These conditions have been 
expressed  in copositivity 
terms as in Simpson and Spector \cite{simspec} for $n=3$, 
Silhavy \cite{silhavy} and Dacorogna \cite{dacorogna} for arbitrary 
$n$ (for an account on the history of results related to this 
problem see the Silhavy or Dacorogna op. cit.). The conditions 
given in theorem \ref{trank1} have some advantages. The relation 
between rank one convexity and Schur convexity, which is rather 
obvious, can be used to obtain quasiconvexity results. As for the 
condition (b), it contains one inequality instead  a $2^{n}$ 
family  of (equally complex) inequalities expressing copositivity. 
Moreover, for $n=2$ or $n=3$, it can be used to obtain explicit 
conditions, as in Dacorogna \cite{dacorogna}. These explicit
 conditions (for $n=3$), contained in theorem 5, Dacorogna, op. 
cit., are clearly not independent and have a rather involved form. 
I think that for practical purposes 
it is much easier to think in other terms. The next proposition, 
with a straightforward proof, is relevant. 

\begin{prop}
Let $H$, $G$ be two symmetric $n\times n$ matrices, such that $G$ 
has positive entries. The following statements are equivalent: 
\begin{enumerate}
\item[(a)] for any $x \in R^{n}$ 
$$H_{ij}x_{i}x_{j} + G_{ij}\mid x_{i}\mid \mid x_{j}\mid \geq 0$$
\item[(b)] we have the set inclusion $A(G) \subset A(-H)$, where
$$A(G) = \left\{ x \in R^{n} \ : \ G_{ij}\mid x_{i}\mid \mid x_{j}
\mid \leq 1 \right\}$$
$$A(-H) = \left\{ x \in R^{n} \ : \ - H_{ij} x_{i}x_{j}\leq 1 
\right\}$$
\end{enumerate}
\end{prop}
The matter of finding conditions upon $H$, $G$ such that inclusion 
(b) happens is one of comparing asymptotic and extremal properties 
of the sets $A(G)$ and $A(-H)$. These sets have simple descriptions 
and the algebraic conditions upon $H$, $G$ reflects nothing but the 
geometrical effort to put $A(G)$ inside $A(-H)$. If 
$H$ is positive definite then $A(-H) = R^{n}$. Otherwise the 
inclusion $A(G) \subset A(-H)$  can be expressed in terms of 
eigenvalues and eigenvectors 
of $H$, $G$. I don't pursue this path here, because it is separate 
from the purpose of this note.

\section{Majorisation and quasiconvexity}

The goal of this section is to give a class of multiplicative 
quasiconvex isotropic functions which seem to be complementary 
to the polyconvex isotropic ones. We quote the following 
result of Thompson and Freede \cite{thomfre}, Ball \cite{ball2} 
(for a proof coherent with this paper see Le Dret \cite{ledret}). 

\begin{thm}
Let $g:[0,\infty)^{n} \rightarrow R$ be convex, symmetric and 
nondecreasing in each variable. Define the function $w$ by
$$w:gl(n,R) \rightarrow R \ , \ \ w(F) = g(\sigma(F)).$$
Then $w$ is convex. 
\label{thomfreede}
\end{thm}

The main result of this section is: 

\begin{thm}
Let $g:(0,\infty)^{n} \rightarrow R$ be a continuous symmetric 
function and  $h:R^{n} \rightarrow R$, $h = g \circ \exp$. Consider 
also the function $p: R^{n} \rightarrow R$ 
$$p(\sum_{i=1}^{k} \down{x}_{i}) \ = \ h(x_{k})$$
Suppose that: 
\begin{enumerate}
\item[(a)] $h$ is convex, 
\item[(b)] $p$ is nonincreasing in each argument. 
\end{enumerate}
  Let $\Omega \subset R^{n}$ be bounded, with piecewise 
smooth boundary and $\phi: \overline{\Omega} \rightarrow R$ be any 
Lipschitz function  such that $D\phi(x) \in GL(n,R)^{+}$ a.e. and 
$\phi(x) = x$ on $\partial \Omega$. 
Define the function
$$w: GL(n,R)^{+} \rightarrow R \ , \ \ w(F) = g(\sigma(F))$$
Then for any $F \in GL(n,R)^{+}$ we have: 
\begin{equation}
\int_{\Omega} w(F D\phi(x)) \ \geq \ \mid \Omega \mid w(F)
\label{firstqc}
\end{equation}
\label{texe}
\end{thm}
The notion of multiplicative quasiconvexity is given further. 

\begin{defi}
Let $w: GL(n,R)^{+} \rightarrow R$ be a function and $\Omega = 
(0,1)^{n}$. $w$ is multiplicative quasiconvex if for any $F \in 
GL(n,R)^{+}$ and for any Lipschitz function $u: \Omega 
\rightarrow R$, such that for almost any $x \in \Omega$ 
$\det Du(x) >0$ and $u(x)=x$ on $\partial \Omega$, we have the 
inequality: 
$$\int_{\Omega} w(F Du(x)) \ \geq \ \int_{\Omega}w(F)$$
\label{dmq}
\end{defi}

\begin{rk}
In the above definition $\Omega$ can be replaced by any bounded 
open set with piecewise smooth boundary. 
\end{rk}

The notion of multiplicative quasiconvexity appears with the name 
Diff-quasiconvexity in Giaquinta, Modica, Soucek \cite{gms}, page 
174, definition 3. It can be found for the first time in Ball 
\cite{ball2}, under a disguised form. 
It is in fact the natural notion to be considered in connection 
with continuous media mechanics. Any polyconvex function is multiplicative polyconvex. Moreover, classical quasiconvexity implies multiplicative quasiconvexity. Conversely, multiplicative quasiconvexity means quasiconvexity if one extends $w$ on the whole 
$gl(n,R)$ by $w(F) = + \infty$ if $\det F \leq 0$. For the lower 
semicontinuity properties of multiplicative quasiconvex functions 
see Buliga \cite{buliga}. 

Theorem \ref{texe} tells that any $w$ which satisfies the 
hypothesis is multiplicative quasiconvex.

In order to prepare the proof of theorem \ref{texe}, two lemmas 
are given. 

\begin{lema}
Let $h:R^{n} \rightarrow R$ be continuous,  Schur convex and $g = 
h \circ \log$. Define 
$$w: GL(n,R)^{+} \rightarrow R \ , \ \ w(F) = g(\sigma(F))$$
$$\tilde{w}: GL(n,C) \rightarrow R \ , \ \ \tilde{w}(F) = g(\mid 
\lambda(F) \mid)$$
Then for any $F$ 
$$w(F) \geq \tilde{w}(F)$$
\label{l1}
\end{lema}

\begin{proof}
This is a straightforward consequence of the Weyl inequality 
(theorem \ref{tweyl})
$$\log \mid \lambda(F) \mid \ \prec \ \log \sigma(F)$$
and of the Schur convexity of $h$. 
\end{proof}

\begin{lema}
With the notations from the lemma \ref{l1}, for any two symmetric 
matrices $A,B$, we have 
$$\tilde{w}(\exp A \ \exp B) \geq \tilde{w}(\exp(A+B))$$
\label{l2}
\end{lema}

\begin{proof}
We have to check the conditions from Thompson \cite{thompson}, 
Lemma 6, which gives 
sufficient conditions on the function $\tilde{w}$ in order to 
satisfy the inequality we are trying to prove. These conditions are: 
\begin{enumerate}
\item[(1)] for any $X$ and any symmetric positive definite $Y$  
$\tilde{w}(XY) = \tilde{w}(YX)$. 
This is satisfied by definition of $\tilde{w}$. 
\item[(2)] for any $X$ and any $m = 1, 2, ...$ 
$$\tilde{w}\left(\left[XX^{*}\right]^{m}\right) \geq 
\tilde{w}\left(X^{2m}\right)$$
From the definition of $\tilde{w}$ and Lemma \ref{l1} we find 
that $\tilde{w}$ satisfies this condition too. 
\end{enumerate}
\end{proof}

We give now the proof of the theorem \ref{texe}.
\begin{proof}
To  any $F \in GL(n,R)^{+}$ we associate it's polar decomposition 
$F=R_{F}U_{F} = V_{F}R_{F}$. For any function $\phi$ such that 
$D\phi(x) \in GL(n,R)^{+}$ we shall use the (similar) notation
$$D\phi(x) = R\phi(x) U\phi(x) = V\phi(x) R\phi(x)$$
With the notations from the theorem, we have from the isotropy of 
$w$, hypothesis  (a) and  theorem \ref{t5} that $h$ is Schur 
convex. From  lemma \ref{l1} and lemma \ref{l2} we obtain  the 
chain of inequalities $$\int_{\Omega} w(FD\phi(x)) \ = \ \int_{\Omega} w(U_{F}V\phi(x)) \ \geq \ 
\int_{\Omega} \tilde{w}(U_{F}V\phi(x)) \ \geq \ \int_{\Omega} 
\tilde{w}
\left( \exp \left( \log U_{F} + \log V\phi(x) \right) \right)$$
The chain of inequalities continues by using the convexity 
hypothesis (a) 
(suppose that $\mid \Omega \mid = 1$): 
$$\int_{\Omega} \tilde{w}
\left( \exp \left( \log U_{F} + \log V\phi(x) \right) \right) \ 
\geq \ 
\tilde{w}\left( \exp \left(\log U_{F} + \int_{\Omega}\log V\phi(x)
\right) \right)$$
Now, I claim that the matrix 
$$\int_{\Omega}\log V\phi(x)$$
is negative definite. Then, from theorem 
\ref{tf1}, we find that 
$$\down{\lambda} \left(\log U_{F} + \int_{\Omega}\log 
V\phi(x)\right) \leq 
\down{\lambda} (\log U)$$
We use now the nonincreasing condition (b) to finish the chain of 
inequalities 
$$\tilde{w}\left( \exp \left(\log U_{F} + \int_{\Omega}\log V\phi(x)\right)\right) \ \geq 
\ \tilde{w}(U_{F}) = w(U_{F}) = w(F)$$

Let us see, finally, why the matrix $\int_{\Omega}\log V\phi(x)$ 
is negative definite. The function $\log \down{\lambda_{1}}(F)$ 
is well known polyconcave, hence all the functions $\log 
\down{\lambda_{i}}(F)$ satisfy the inequality: 
$$\int_{\Omega} \log \down{\lambda_{i}}(D\phi(x)) \ \leq \ \mid 
\Omega \mid 
\log \down{\lambda_{i}} (I_{n}) = 0$$
Take now any vector $v \in R^{n}$, $v \not = 0$. Remember that 
$V\phi(x)$ is a symmetric matrix which admits the decomposition 
$$V\phi(x) = Q\phi(x) \ Diag (\lambda(D\phi(x))) Q^{T} \phi(x)$$ 
hence 
$$\log V\phi(x) \ = \ Q\phi(x) \ Diag ( \log \lambda(D\phi(x))) 
Q^{T}\phi(x)$$
Therefore 
$$\sum_{i,j=1}^{n}[ \log V\phi(x) ]_{ij} v_{i}v_{j} \ \leq \ \log 
\down{\lambda_{1}} \mid v \mid^{2}$$
Use the inequality given by polyconcavity to deduce the claim. 
\end{proof}

A consequence of the theorem \ref{texe} is:

\begin{prop}
In the hypothesis of the theorem \ref{texe}, the function $w$ is 
rank one convex. 
\label{prank}
\end{prop}

The family of functions satisfying the hypothesis of theorem 
\ref{texe} is non void. Two examples  are given further. 

For the first example take the polar decomposition $F= R_{F}U_{F}$ 
and define the function: $w(F) = \ \log \ trace \ U_{F}^{-1}$. It 
satisfies the hypothesis, by straightforward computation. Indeed, 
using the  notations of theorem 
\ref{texe}, the associated function $g: (0,+\infty)^{n} \rightarrow 
R$ is 
$$g(y_{1},...,y_{n}) \ = \ \log \left( \sum_{i=1}^{n} 
\frac{1}{y_{i}} \right)$$
hence the function $h(x) = g(\exp x)$ has the expression: 
$$h(x_{1},...,x_{n}) \ = \ \log \left( \sum_{i=1}^{n} \exp 
(-x_{i}) \right)$$
which is easy to check that is convex and the associated
function $p$ is decreasing in each argument. 

For the second example consider a modified Ogden potential. Set 
$$\| F \|_{k} \ = \ \left( \prod_{i=1}^{k} \down(\sigma)_{i}(F)
\right)^{1/k}$$ and define: 
$$w(F) \ = \ \sum_{i= 1}^{n} \frac{1}{\| F \|_{i}^{\alpha}}$$
for some $\alpha \geq 2$. The associated function $h$ is then 
$$h(x_{1}, ... , x_{n}) \ = \ \sum_{k = 1}^{n} \exp \left(
(-\alpha/k) \sum_{i=1}^{k} \down{x}_{i} \right)$$
which again satisfies the hypothesis of the theorem.

Both functions are not known to be polyconvex. In fact they do not 
satisfy the following sufficient condition for polyconvexity, due 
to Ball \cite{ball2}, given here for simplicity for $n=3$ (see also 
Le Dret \cite{ledret}). 

\begin{thm}
Let $\phi: R^{3}_{+} \times R^{3}_{+} \times R_{+} \rightarrow R$ 
which is nondecreasing in the first six variables and such that for 
any pair of permutations $\sigma,\tau \in S_{3}$ 
$$\phi(v_{\sigma(1)}, v_{\sigma(2)}, v_{\sigma(3)}, v_{\tau(1)+3},
v_{\tau(2)+3},v_{\tau(3)+3}, v_{7}) \ = \ \phi(v_{1},...,v_{7})$$
Then the function: 
$$w(F) \ = \ \phi(\sigma_{1}(F),...,
\sigma_{1}(F)\sigma_{2}(F),..., 
\det F)$$
is polyconvex. 
\label{critball}
\end{thm}

This is the reason for thinking that the functions described in 
theorem \ref{texe} are complementary to objective, isotropic, 
polyconvex ones. Remark nevertheless 
that the function 
$$w(F) \ = \ - \log \det F$$
satisfies the hypothesis of theorem \ref{texe} and it is also 
polyconvex. 

Let us consider only the Schur convexity and componentwise 
convexity hypothesis related to $w$. 

\begin{prop}
Let $h: R^{n} \rightarrow R$ be Schur convex and the function 
$x \in R \mapsto h(\log (x), ..., \log(x))$ be convex, continuous. 
Let $\phi: \Omega \rightarrow R$ be such that 
almost everywhere we have $D\phi(x) \in GL(n,R)^{+}$, 
$$\int_{\Omega} D\phi(x) \ = I_{n}$$
and the map $x \mapsto w(D\phi(x))$ is integrable. 
Then 
$$\int_{\Omega} w(D\phi(x)) \ \geq \ \mid \Omega \mid w(I_{n})$$
\label{pschur}
\end{prop}

 \begin{proof}
Because $h$ is Schur convex and for almost any $x \in \Omega$  
$$\frac{1}{n}\log \det D\phi(x) (1,...1) \prec \ \log 
\sigma(D\phi(x))$$
we have the inequality 
$$w(D\phi(x)) \geq w\left( (\det D\phi(x))^{1/n} I_{n}\right)$$
Use the convexity hypothesis to obtain the desired inequality. 
\end{proof}

\section{Rank one convex hulls and majorisation}

In this section it is explained how majorisation appears in the 
representation of some rank one convex hulls. What would be really 
nice to understand are the implications of Lie group aspects of 
majorisation onto calculus of variations. 
The fact that such implications should exist is straightforward, 
but far from being self evident. 
 
Further is given a  proof of theorem \ref{dacotan} using 
majorisation. In
 this proof we use the fact that majorisation relation 
$$x \prec \prec y \ \mbox{ if } \ \log x \prec \log y$$
 is defined using polyconvex maps. The isotropy of the set 
$E(a)$ from theorem \ref{dacotan} implies that the description of 
it's rank one convex hull reduces to the description of the set of 
matrices $B \prec Diag(a)$, where $\prec$ is Thompson's order 
relation. These facts (partially) explain 
the resemblance between theorems \ref{hornthom} and \ref{dacotan}. 

Let $a \in (0,\infty)^{n}$. Denote by $E(a)$ the set of matrices 
$F$ with positive determinant 
such that $\sigma(F) = Pa$ for some $P \in S_{n}$. We have to prove 
that 
$$Pco \ E(a) \ = \ Rco \ E(a) \ = \ K(a)$$
where 
$$K(a) \ = \ \left\{ B \in GL(n,R)^{+}  \ : \ B \prec Diag \ (a) 
\right\}$$
The set $K(a)$ is polyconvex, being an intersection of preimages of 
$(-\infty, 0]$ by polyconvex functions. Therefore 
$$Rco \ E(a) \ \subset \ Pco \ E(a) \ \subset \ K(a)$$
It is left to prove that $K(a) \subset  Rco \ E(a)$. For this 
remark that $E(a)$ can be written as: 
$$E(a) \ = \ \left\{ R \ P.Diag(a) \ Q \ : \ R,Q \in SO(n) \ , \ P 
\in S_{n} \right\}$$
Consider the convex cone of functions ($Rco$ denotes the class of 
rank one convex functions)
$$Rco (a) \ = \ \left\{ \phi \in Rco \ : \ \forall A \in E(a) \ \ 
\phi(A) = 0 
\right\}$$
This cone is closed to $sup$ operation. Moreover, it has the same 
symmetries as 
$E(a)$, that is for any $R,Q \in SO(n)$ and any $P \in S_{n}$ we have 
$$\phi \in \ Rco (a) \ \Longrightarrow \ \left[ F \in GL(n,R)^{+} 
\mapsto (R,Q,P).\phi(F) = \phi(R \ P.F \ Q) \right] \in Rco(a)$$
Hence if $\phi \in Rco(a)$ then $\bar{\phi} \in Rco(a)$, where 
$\bar{\phi}$ is the objective isotropic function 
$$\bar{\phi}(F) \ = \ \sup \left\{ (R,Q,P).\phi(F) \ : \ R,Q \in 
SO(n) \ , \ P \in S_{n} \right\}$$
Objective isotropic rank one convex functions have Schur convex 
diagonal, as a consequence of theorem \ref{trank1} (a) (if the rank 
one convex $w$ is not $C^{2}$ use a convolution argument).  
Therefore $F \in K(a)$ and 
$\phi \in Rco(a)$ imply 
$$\phi(F) \leq \bar{\phi}(F) \leq \bar{\phi}(Diag(a)) = 0$$
This proves the inclusion $K(a) \subset Rco(a)$.
 
{\bf Aknowledgements.} I have learned about Horn-Thompson theorem 
with its amazing implications in Lie group convexity results from 
discussions with Tudor Ratiu. I want to thank Bernard Dacorogna for 
keeping me connected with parts of his research.

\vspace{2.cm}

\noindent
Institute of Mathematics of the Romanian Academy \\ 
 PO BOX 1-764, RO 70700, Bucharest, Romania, 
 e-mail: Marius.Buliga@imar.ro \\ and \\ 
Ecole Polytechnique F\'ed\'erale de Lausanne, D\'epartement de 
Math\'ematiques \\ 
1015 Lausanne, Switzerland, e-mail: Marius.Buliga@epfl.ch

\end{document}